\documentclass[12pt,draft]{amsart}

\usepackage[cp1251]{inputenc}
\usepackage[T2A]{fontenc}
\usepackage[english]{babel}
\usepackage{amsmath,amsfonts,amssymb}
\usepackage{geometry}
\newtheorem{theorem}{Theorem}[section]

\newtheorem{proposition}{Proposition}[section]

\theoremstyle{remark}

\theoremstyle{definition}
\newtheorem{condition}{Condition}

\title[Elliptic problems and H\"ormander spaces]
{Elliptic problems and H\"{o}rmander spaces}

\author[V. A. Mikhailets]{Vladimir A. Mikhailets}

\address{%
Institute of Mathematics, National Academy of Sciences of Ukraine\\
Tereshchenkivs'ka str. 3\\
01601 Kyiv\\
Ukraine}

\email{mikhailets@imath.kiev.ua}

\author[A. A. Murach]{Aleksandr A. Murach}

\address{%
Institute of Mathematics, National Academy of Sciences of Ukraine\\
Tereshchenkivs'ka str. 3\\
01601 Kyiv\\
Ukraine\\
Chernigiv State Technological University\\
Shevchenka str. 95\\
14027 Chernigiv\\
Ukraine}

\email{murach@imath.kiev.ua}

\subjclass[2000]{Primary 35J30, 35J40, Secondary 46E35}

\keywords{H\"ormander spaces, generalized smoothness, interpolation with a
function parameter, elliptic operator, elliptic boundary-value problem, the
Fredholm property, local regularity of solutions, the Lions-Magenes theorems}

\begin{document}

\maketitle

\begin{abstract}
The paper gives a survey of the modern results on elliptic problems on the
H\"ormander function spaces. More precisely, elliptic problems are studied on a
Hilbert scale of the isotropic H\"ormander spaces parametrized by a real number
and a function slowly varying at $+\infty$ in the Karamata sense. This refined
scale is finer than the Sobolev scale and is closed with respect to the
interpolation with a function parameter. The Fredholm property of elliptic
operators and elliptic boundary-value problems is preserved for this scale.
A~local refined smoothness of the elliptic problem solution is studied. An
abstract construction of classes of function spaces in which the elliptic
problem is a Fredholm one is found. In particular, some generalizations of the
Lions-Magenes theorems are given.
\end{abstract}

\section*{0. Introduction}

The paper gives a survey of the modern results [32--49] devoted to elliptic
problems on the Hilbert scale of the isotropic H\"ormander spaces
$$
H^{s,\varphi}:=H_{2}^{\langle\cdot\rangle^{s}\,\varphi(\langle\cdot\rangle)},
\quad\langle\xi\rangle:=\bigl(1+|\xi|^{2}\bigr)^{1/2}. \eqno(0.1)
$$
Here $s\in\mathbb{R}$ and $\varphi$ is a functional parameter slowly varying at
$+\infty$ in the Karamata sense. In particular, every standard function
$$
\varphi(t)=(\log t)^{r_{1}}(\log\log t)^{r_{2}}\ldots(\log\ldots\log
t)^{r_{k}},\quad t\gg1,
$$
$$
\{r_{1},r_{2},\ldots,r_{k}\}\subset\mathbb{R},\;k\in\mathbb{Z}_{+},
$$
is admissible. This scale contains the Sobolev scale
$\{H^{s}\}\equiv\{H^{s,1}\}$, is attached to it by the number parameter $s$,
and much finer than $\{H^{s}\}$.

Spaces of form (1) arise naturally in different spectral problems: convergence
of spectral expansions of self-adjoint elliptic operators almost everywhere, in
the norm of the spaces $L_{p}$ with $p>2$ or $C$ (see survey \cite{AIN76});
spectral asymptotics of general self-adjoint elliptic operators in a bounded
domain, the Weyl formula, a sharp estimate of the remainder in it (see
\cite{Mikhailets82, Mikhailets89}) and others. They may be expected to be
useful in other "fine"\, questions.  Due to their interpolation properties, the
spaces $H^{s,\varphi}$ occupy a special position among the spaces of a
generalized smoothness, which are actively investigated and used today (see
survey \cite{KalLiz87}, recent articles \cite{HarMou04, FarLeo06} and the
bibliography given therein).

The paper consists of six sections. In Section~1 the refined scale of the
H\"ormander spaces (0.1) is introduced and studied. In particular, important
interpolation properties of this scale are under investigation. In Section~2 an
elliptic pseudodifferential operator on the refined scale on a closed compact
smooth manifold is considered. We show that this operator is a Fredholm one and
establishes a collection of isomorphisms on the two-sided refined scale. The
local refined smoothness of a solution to the elliptic equation is studied. We
also give an equivalent definition of the refined scale on the closed manifold
by means of certain functions of a positive elliptic operator.

Next we study a regular elliptic boundary problem on a bounded Euclidean domain
with the smooth boundary. In Section~3 we show that the operator of this
problem is a Fredholm one on the upper part of the refined scale. A local
refined smoothness up to the boundary of a solution to the problem is studied.
As an important application, we give a sufficient condition for the solution to
be classical. Section~4 is devoted to semihomogeneous elliptic boundary
problems. We show that these problems are Fredholm on the two-sided refined
scales.

Since the operator of the general nonhomogeneous boundary problem cannot be
defined correctly on the lower part of the refined scale, we consider in
Section~5 a special modified refined scale on which the operator is
well-defined, bounded, and Fredholm everywhere. This modification depends
solely on the order of the problem, so that the theorem on the Fredholm
property is generic for the class of elliptic problems having the same order.

The last Section~6 is devoted to some individual theorems on the Fredholm
property. We give an abstract construction of classes of function spaces on
which the elliptic problem operator is a Fredholm one. A characteristic feature
of this construction is that the domain of the operator depends on coefficients
of the elliptic expression. So, we have the individual theorems on the Fredholm
property. As an important application, we give some generalizations of the
known Lions-Magenes theorems.

\section{A refined scale of H\"ormander spaces} 

Let us denote by $\mathcal{M}$ the set of all functions
$\varphi:[1,+\infty)\rightarrow(0,+\infty)$ such that:
\begin{itemize}
\item [a)] $\varphi$ is a Borel measurable function;
\item [b)] the functions $\varphi$ and $1/\varphi$ are bounded on every closed
interval $[1,b]$, where $1<b<+\infty$;
\item [c)] $\varphi$ is a slowly varying function at $+\infty$
in the Karamata sense (see \cite[Sec. 1.1]{Seneta76}), i.e.
$$
\lim_{t\rightarrow\,+\infty}{\varphi(\lambda\,t)}/{\varphi(t)}=1 \quad\mbox{for
each}\quad\lambda>0.
$$
\end{itemize}

Let $s\in\mathbb{R}$ and $\varphi\in\mathcal{M}$. We denote by
$H^{s,\varphi}(\mathbb{R}^{n})$ the space of all tempered distributions $w$ on
the Euclidean space $\mathbb{R}^{n}$ such that the Fourier transform
$\widehat{w}$ of the distribution $w$ is a locally Lebesgue integrable on
$\mathbb{R}^{n}$ function which satisfies the condition
$$
\int_{\mathbb{R}^{n}}\langle\xi\rangle^{2s}\,\varphi^{2}(\langle\xi\rangle)
\:|\widehat{w}(\xi)|^{2}\,d\xi<\infty.
$$
Here $\langle\xi\rangle=(1+\xi_{1}^{2}+\ldots+\xi_{n}^{2})^{1/2}$ is the
smoothed modulus of a vector $\xi=(\xi_{1},\ldots,\xi_{n})\in \mathbb{R}^{n}$.
An inner product in the space $\mathrm{H}^{s,\varphi}(\mathbb{R}^{n})$ is
defined by the formula
$$
(w_{1},w_{2})_{\mathrm{H}^{s,\varphi}(\mathbb{R}^{n})}:=
\int_{\mathbb{R}^{n}}\langle\xi\rangle^{2s}\varphi^{2}(\langle\xi\rangle)
\,\widehat{w_{1}}(\xi)\,\overline{\widehat{w_{2}}(\xi)}\,d\xi.
$$
The inner product induces the norm in $\mathrm{H}^{s,\varphi}(\mathbb{R}^{n})$
in the usual way. Note that we consider distributions which are antilinear
functionals on the space of test functions.

The space $H^{s,\varphi}(\mathbb{R}^{n})$ is a special isotropic Hilbert case
of the spaces introduced and investigated by L.~H\"ormander
\cite[Sec.~2.2]{Hermander63}, \cite[Sec.~10.1]{Hermander83} and the different
spaces studied by L.~R.~Volevich and B.~P.~Paneah \cite[Sec.~2]{VoPa65},
\cite[Sec.~1.4.2]{Paneah00}. In the simplest case where
$\varphi(\cdot)\equiv1$, the space $H^{s,\varphi}(\mathbb{R}^{n})$ coincides
with the Sobolev space $H^{s}(\mathbb{R}^{n})$. The inclusions
$$
\bigcup_{\varepsilon>0}H^{s+\varepsilon}(\mathbb{R}^{n})=:H^{s+}(\mathbb{R}^{n})
\subset H^{s,\varphi}(\mathbb{R}^{n})\subset
H^{s-}(\mathbb{R}^{n}):=\bigcap_{\varepsilon>0}H^{s-\varepsilon}(\mathbb{R}^{n})
$$
imply that in the set of separable Hilbert spaces
$$
\bigl\{H^{s,\varphi}(\mathbb{R}^{n}):s\in\mathbb{R},\varphi\in\mathcal{M}\,\bigr\},
\eqno(1.1)
$$
the functional parameter $\varphi$ defines an additional (subpower) smoothness
with respect to the basic (power) $s$-smoothness. Otherwise speaking, $\varphi$
\textit{refines} the power smoothness $s$. Therefore, the collection of spaces
(1.1) is naturally called the \textit{refined} scale over $\mathbb{R}^{n}$
(with respect to the Sobolev scale).

We are going to study an application of the refined scale to elliptic boundary
problems in a bounded domain $\Omega\subset\mathbb{R}^{n}$. Therefore, we need
to have the refined scales over the domain $\Omega$ and over its boundary
$\partial\Omega$. The refined scale over the closed domain
$\overline{\Omega}:=\Omega\cup\partial\Omega$ is also of use. We construct
these scales from (1.1) in the standard way.

Let us denote
$$
H^{s,\varphi}(\Omega):=\bigl\{u=w\upharpoonright\Omega:\,w\in
H^{s,\varphi}(\mathbb{R}^{n})\bigr\},
$$
$$
\|\,u\,\|_{H^{s,\varphi}(\Omega)}:=
\inf\,\bigl\{\,\|\,w\,\|_{H^{s,\varphi}(\mathbb{R}^{n})}:\,w\in
H^{s,\varphi}(\mathbb{R}^{n}),\;\;w=u\;\;\mbox{in}\;\;\Omega\,\bigr\}.
$$
The norm in the space $H^{s,\varphi}(\Omega)$ is induced by the inner product
$$
\bigl(u_{1},u_{2}\bigr)_{H^{s,\varphi}(\Omega)}:= \bigl(w_{1}-\Pi
w_{1},w_{2}-\Pi w_{2}\bigr)_{H^{s,\varphi}(\mathbb{R}^{n})}.
$$
Here $w_{j}\in H^{s,\varphi}(\mathbb{R}^{n})$, $w_{j}=u_{j}$ in $\Omega$ for
$j=1,\,2$, and $\Pi$ is the orthogonal projector of the space
$H^{s,\varphi}(\mathbb{R}^{n})$ onto the subspace $\{w\in
H^{s,\varphi}(\mathbb{R}^{n}):\mathrm{supp}\,w\subseteq\mathbb{R}^{n}
\setminus\Omega\}$. The space $H^{s,\varphi}(\Omega)$ is a separable Hilbert
one.

We also denote
$$
H^{s,\varphi}_{\overline{\Omega}}(\mathbb{R}^{n}):=\bigl\{w\in
H^{s,\varphi}(\mathbb{R}^{n}):\,\mathrm{supp}\,w\subseteq
\overline{\Omega}\,\bigr\}.
$$
This space is a separable Hilbert one with respect to the inner product in the
space $H^{s,\varphi}(\mathbb{R}^{n})$.

Thus the space $H^{s,\varphi}(\Omega)$ consists of the distributions given in
the open domain $\Omega$, whereas the space
$H^{s,\varphi}_{\overline{\Omega}}(\mathbb{R}^{n})$ consists of the
distributions supported on the closed domain $\overline{\Omega}$. The
collections of Hilbert spaces
$$
\bigl\{H^{s,\varphi}(\Omega):s\in\mathbb{R},\varphi\in\mathcal{M}\,\bigr\}
\quad\mbox{and}\quad \bigl\{H^{s,\varphi}_{\overline{\Omega}}(\mathbb{R}^{n}):
s\in\mathbb{R},\varphi\in\mathcal{M}\,\bigr\} \eqno(1.2)
$$
are called the refined scales over $\Omega$ and over $\overline{\Omega}$
respectively.

The boundary $\partial\Omega$ is assumed to possess an infinitely smooth field
of unit vectors of normals. So, $\partial\Omega$ is a particular case of a
compact closed infinitely smooth manifold. Let us define the refined scale over
a closed infinitely smooth manifold $\Gamma$ of an arbitrary dimension $n$.

We choose a finite atlas from the $C^{\infty}$-structure on the manifold
$\Gamma$ consisting of the local charts
$\alpha_{j}:\mathbb{R}^{n}\leftrightarrow U_{j}$, $j=1,\ldots,r$. Here the open
sets $U_{j}$ form the finite covering of the manifold $\Gamma$. Let functions
$\chi_{j}\in C^{\infty}(\Gamma)$, $j=1,\ldots,r$, form a partition of unity on
$\Gamma$ satisfying the condition $\mathrm{supp}\,\chi_{j}\subset U_{j}$.

We set
$$
H^{s,\varphi}(\Gamma):=\left\{h\in\mathcal{D}'(\Gamma):\;
(\chi_{j}h)\circ\alpha_{j}\in
H^{s,\varphi}(\mathbb{R}^{n})\;\;\forall\;j=1,\ldots,r\right\}.
$$
Here, as usual, $\mathcal{D}'(\Gamma)$ is the topological space of all
distributions on $\Gamma$, and $(\chi_{j}h)\circ\alpha_{j}$ is the
representation of the distribution $\chi_{j}h$ in the local chart $\alpha_{j}$.
The inner product in the space $H^{s,\varphi}(\Gamma)$ is defined by the
formula
$$
(h_{1},h_{2})_{H^{s,\varphi}(\Gamma)}:=
\sum_{j=1}^{r}\,((\chi_{j}h_{1})\circ\alpha_{j},
(\chi_{j}h_{2})\circ\alpha_{j})_{H^{s,\varphi}(\mathbb{R}^{n})}
$$
and induces the norm in the usual way.

The Hilbert space $H^{s,\varphi}(\Gamma)$ is separable and does not depend (up
to equivalence of norms) on the choice of the atlas and the partition of unity.
The collection of function spaces
$$
\{H^{s,\varphi}(\Gamma):s\in\mathbb{R},\;\varphi\in\mathcal{M}\} \eqno(1.3)
$$
is called the refined scale over the manifold $\Gamma$. Specifically, we need
the refined scale of spaces $H^{s,\varphi}(\partial\Omega)$.

We note the following properties of the refined scales:

\begin{theorem} 
Let $s\in\mathbb{R}$ and $\varphi,\varphi_{1}\in\mathcal{M}$. The following
assertions are true:
\begin{itemize}
\item [(i)] The set $C^{\infty}(\,\overline{\Omega}\,)$ is dense in the space
$H^{s,\varphi}(\Omega)$.
\item [(ii)] The set $C^{\infty}_{0}(\Omega):=\{w\in C^{\infty}(\mathbb{R}^{n}):
\mathrm{supp}\,w\subset\Omega\}$ is dense in the space
$H^{s,\varphi}_{\overline{\Omega}}(\mathbb{R}^{n})$.
\item [(iii)] If $|s|<1/2$, then the mapping $w\rightarrow w\upharpoonright\Omega$
establishes a topological isomorphism from
$H^{s,\varphi}_{\overline{\Omega}}(\mathbb{R}^{n})$ onto
$H^{s,\varphi}(\Omega)$.
\item [(iv)] For each $\varepsilon>0$ the compact and dense embeddings hold:
$$
H^{s+\varepsilon,\,\varphi_{1}}(\Omega)\hookrightarrow H^{s,\varphi}(\Omega),
\quad
H^{s+\varepsilon,\,\varphi_{1}}_{\overline{\Omega}}(\mathbb{R}^{n})\hookrightarrow
H^{s,\varphi}_{\overline{\Omega}}(\mathbb{R}^{n}). \eqno(1.4)
$$
\item [(v)] Suppose that the function $\varphi/\varphi_{1}$
is bounded in a neighborhood of $+\infty$. Then continuous dense embeddings
$(1.4)$ are valid for $\varepsilon=0$. They are compact if
$\varphi(t)/\varphi_{1}(t)\rightarrow0$ as $t\rightarrow+\infty$.
\item [(vi)] For every fixed integer $k\geq0$ the inequality
$$
\int_{1}^{\,+\infty}\frac{d\,t}{t\,\varphi^{\,2}(t)}<\infty \eqno(1.5)
$$
is equivalent to the embedding $H^{k+n/2,\,\varphi}(\Omega)\hookrightarrow
C^{k}(\,\overline{\Omega}\,)$. This embedding is compact.
\item [(vii)] The spaces $H^{s,\varphi}(\Omega)$ and
$H^{-s,1/\varphi}_{\overline{\Omega}}(\mathbb{R}^{n})$ are mutually dual with
respect to the inner product in $L_{2}(\Omega)$.
\item [(viii)] The mapping $u\mapsto u\upharpoonright\partial\Omega$,
$u\in C^{\infty}(\,\overline{\Omega}\,)$, is extended by a continuity to the
bounded trace operator from $H^{s,\varphi}(\Omega)$ onto
$H^{s-1/2,\,\varphi}(\partial\Omega)$, provided that $s>1/2$.
\end{itemize}
\end{theorem}

Assertions (iv) -- (vi) show that the refined scale is much finer than the
classical Sobolev scale (the case of $\varphi\equiv\varphi_{1}\equiv1$). Note
also that $\varphi\in\mathcal{M}\Leftrightarrow1/\varphi\in\mathcal{M}$, so the
space $H^{-s,1/\varphi}_{\overline{\Omega}}(\mathbb{R}^{n})$ in assertion (vii)
is defined as an element of the refined scale.

\begin{theorem} 
Let $s\in\mathbb{R}$ and $\varphi,\varphi_{1}\in\mathcal{M}$. Then:
\begin{itemize}
\item [(i)] Assertions \rm (i) \it and \rm (iv) \it -- \rm (vi) \it of
Theorem $1.1$ hold true if we replace both the notations $(\Omega)$ and
$(\,\overline{\Omega}\,)$ with $(\Gamma)$.
\item [(ii)] The spaces $H^{s,\varphi}(\Gamma)$ and
$H^{-s,1/\varphi}(\Gamma)$ are mutually dual (up to equivalence of norms) with
respect to the inner product in the space $L_{2}(\Gamma,dx)$, where $dx$ is a
$C^{\infty}$-smooth density on $\Gamma$.
\end{itemize}
\end{theorem}

The refined scale of spaces (1.1), (1.2), and (1.3) were introduced and
investigated by authors in \cite{MM1, MM2b, MM6}. Theorems 1.1, 1.2 were proved
in \cite[Theorem 3.6]{MM2b} and \cite[Theorem 4.2]{MM6}. All assertions of
these theorems, except (iii), follow from the properties of H\"ormander spaces
\cite[Sec.~2.2]{Hermander63}, \cite[Sec.~10.1]{Hermander83} (see also
\cite[Sec.~2]{VoPa65}, \cite[Sec.~1.4.2]{Paneah00}).

The refined scale possesses the interpolation property which selects the scale
from among the spaces of generalized smoothness. Namely, every space of this
scale is obtained by the interpolation, with an appropriate function parameter,
of a couple of the Sobolev spaces. We recall the definition of such an
interpolation in the case of general separable Hilbert spaces.

Let an ordered couple $X:=[X_{0},X_{1}]$ of complex Hilbert spaces $X_{0}$ and
$X_{1}$ be such that these spaces are separable and the continuous dense
embedding $X_{1}\hookrightarrow X_{0}$ holds true. We call this couple
admissible. For the couple $X$ there exists an isometric isomorphism
$J:X_{1}\leftrightarrow X_{\,0}$ such that $J$ is a self-adjoint positive
operator in the space $X_{\,0}$ with the domain $X_{1}$. This operator is
uniquely determined by the couple $X$. Let a Borel measurable function
$\psi:(0,+\infty)\rightarrow (0,+\infty)$ be given. We denote by
$[X_{0},X_{1}]_{\psi}$ or simply by $X_{\psi}$ the domain of the operator
$\psi(J)$ endowed with the graphics inner product and the corresponding norm:
$$
(u,v)_{X_{\psi}}:=(u,v)_{X_{0}}+(\psi(J)u,\psi(J)v)_{X_{0}},\quad
\|\,u\,\|_{X_{\psi}}=(u,u)_{X_{\psi}}^{1/2}.
$$
The space $X_{\psi}$ is a separable Hilbert one.

The function $\psi$ is called an \textit{interpolation parameter} if the
following condition is fulfilled for all admissible couples $X=[X_{0},X_{1}]$,
$Y=[Y_{0},Y_{1}]$ of Hilbert spaces and an arbitrary linear mapping $T$ given
on $X_{0}$: if the restriction of the mapping $T$ to the space $X_{j}$ is a
bounded operator $T:X_{j}\rightarrow Y_{j}$ for each $j=0,\,1$, then the
restriction of the mapping $T$ to the space $X_{\psi}$ is also a bounded
operator $T:X_{\psi}\rightarrow Y_{\psi}$.

\begin{theorem} 
Let a function $\varphi\in\mathcal{M}$ and positive numbers
$\varepsilon,\delta$ be given. We set
$$
\psi(t):=t^{\,\varepsilon/(\varepsilon+\delta)}\,\varphi(t^{1/(\varepsilon+\delta)})
\;\;\mbox{for}\;\;t\geq1\quad\mbox{and}\quad\psi(t):=\varphi(1)\;\;\mbox{for}\;\;
0<t<1.
$$
Then the function $\psi$ is an interpolation parameter and, for each
$s\in\mathbb{R}$, the following equalities of spaces with equivalence of norms
in them are true:
$$
\bigl[H^{s-\varepsilon,1}(G),H^{s+\delta,1}(G)\bigr]_{\psi}=H^{s,\varphi}(G)\quad
\mbox{for}\quad G\in\{\mathbb{R}^{n},\Omega,\Gamma\},
$$
$$
\bigl[H^{s-\varepsilon,1}_{\overline{\Omega}}(\mathbb{R}^{n}),
H^{s+\delta,1}_{\overline{\Omega}}(\mathbb{R}^{n})\bigr]_{\psi}=
H^{s,\varphi}_{\overline{\Omega}}(\mathbb{R}^{n}).
$$
\end{theorem}

The refined scale is closed with respect to the interpolation with a function
parameter $\psi(t):=t^{\theta}\chi(t)$ where $0<\theta<1$, whereas $\chi(t)$ is
a Borel measurable positive function slowly varying at $+\infty$.

\begin{theorem} 
Let $s_{0},s_{1}\in\mathbb{R}$, $s_{0}\leq s_{1}$, and
$\varphi_{0},\varphi_{1}\in\mathcal{M}$. In the case where $s_{0}=s_{1}$ we
suppose that the function $\varphi_{0}/\varphi_{1}$ is bounded in a
neighborhood of $+\infty$. Let a Borel measurable function
$\psi:(0,+\infty)\rightarrow(0,+\infty)$ is of the form
$\psi(t):=t^{\theta}\chi(t)$, where $0<\theta<1$ and $\chi(t)$ is a function
slowly varying at $+\infty$. Then $\psi$ is an interpolation parameter, and the
following equalities of spaces with equivalence of norms in them are true:
$$
\bigl[H^{s_{0},\varphi_{0}}(G), H^{s_{1},\varphi_{1}}(G)\bigr]_{\psi}=
H^{s,\varphi}(G)\quad \mbox{for}\quad G\in\{\mathbb{R}^{n},\Omega,\Gamma\},
$$
$$
\bigl[H^{s_{0},\varphi_{0}}_{\overline{\Omega}}(\mathbb{R}^{n}),
H^{s_{1},\varphi_{1}}_{\overline{\Omega}}(\mathbb{R}^{n})\bigr]_{\psi}=
H^{s,\varphi}_{\overline{\Omega}}(\mathbb{R}^{n}).
$$
Here $s:=(1-\theta)s_{0}+\theta s_{1}$, and the function
$\varphi\in\mathcal{M}$ is given by the formula
$$
\varphi(t):=\varphi_{0}^{1-\theta}(t)\,\varphi_{1}^{\theta}(t)\,
\chi\left(t^{s_{1}-s_{0}}\varphi_{1}(t)/\varphi_{0}(t)\right)\quad\mbox{for}\quad
t\geq1.
$$
\end{theorem}

The interpolation of general Hilbert spaces with a function parameter was
studied in \cite{FoiasLions61, Donoghue67, Pustylnik82, MM11}. The class of all
interpolation parameters was described in \cite{Pustylnik82} (see also
\cite[Theorem 2.7]{MM11}). Theorem 1.3 was proved in \cite[Theorems 3.1,
3.5]{MM2b} and \cite[Theorem 4.1]{MM6}. Theorem 1.4 was proved in \cite[Theorem
3.7]{MM11} for the refined scale over $\Gamma$ (the proof for the scales (1.1)
and (1.2) is analogous). Various normed spaces of generalized smoothness over
$\mathbb{R}^{n}$ were studied by means of the interpolation with a function
parameter in \cite{Merucci84, CobFern86}.

\section{An elliptic operator on a closed manifold} 

We recall that $\Gamma$ is a closed (compact and without a boundary) infinitely
smooth manifold of an arbitrary dimension $n\geq1$ and a certain
$C^{\infty}$-density $dx$ is defined on $\Gamma$. We interpret
$\mathcal{D}'(\Gamma)$ as a space antidual to $C^{\infty}(\Gamma)$ with respect
to the extension of the inner product in $L_{2}(\Gamma,\mathrm{d}x)$ by
continuity. This extension is denoted by $(f,w)_{\Gamma}$ for
$f\in\mathcal{D}'(\Gamma)$, $w\in C^{\infty}(\Gamma)$.

Let $A$ be a classical (polyhomogeneous) pseudodifferential operator on
$\Gamma$ of an arbitrary order $r\in\mathbb{R}$. The complete symbol of $A$ is
an infinitely smooth complex-valued function on the cotangent bundle
$T^{\ast}\Gamma$. We assume that pseudodifferential operator $A$ is elliptic on
$\Gamma$.

The mapping $u\mapsto Au$ is a linear continuous operator on the space
$\mathcal{D}'(\Gamma)$. We will investigate the restriction of this operator to
spaces of the refined scale over $\Gamma$.

Let us denote by $A^{+}$ a pseudodifferential operator formally adjoint to $A$
with respect to the sesquilinear form $(\cdot,\cdot)_{\Gamma}$. Since both $A$
and $A^{+}$ are elliptic on $\Gamma$, both the spaces
$$
N:=\left\{\,u\in C^{\infty}(\Gamma):
\,Au=0\;\;\mbox{on}\;\;\Gamma\,\right\},\quad N^{+}:=\left\{v\in
C^{\infty}(\Gamma):\,A^{+}v=0 \;\;\mbox{on}\;\;\Gamma\,\right\}
$$
are finite-dimensional.

Let us recall the following: a linear bounded operator $T:X\rightarrow Y$ is
called a \textit{Fredholm} one if its kernel is finite-dimensional and its
range $T(X)$ is closed in the space $Y$ and has the finite codimension therein.
Here $X$ and $Y$ are Hilbert spaces. The Fredholm operator $T$ has the finite
\textit{index} $\mathrm{ind}\,T:=\dim\ker T-\dim(Y/\,T(X))$.

\begin{theorem} 
A restriction of the mapping $u\mapsto Au$, $u\in\mathcal{D}'(\Gamma)$,
establishes the linear bounded operator
$$
A:\,H^{s,\varphi}(\Gamma)\rightarrow H^{s-r,\,\varphi}(\Gamma)\quad\mbox{for
each}\;\;s\in\mathbb{R},\;\varphi\in\mathcal{M}. \eqno(2.1)
$$
This operator is a Fredholm one, has the kernel $N$ and the range
$$
\bigl\{f\in H^{s-r,\,\varphi}(\Gamma):\,(f,v)_{\Gamma}=0\;\;\forall\;\; v\in
N^{+}\bigr\}.
$$
The index of the operator $(2.1)$ is equal to $\dim N-\dim N^{+}$ and does not
depend on $s$ and $\varphi$.
\end{theorem}

\begin{theorem} 
For arbitrarily chosen parameters $s\in\mathbb{R}$, $\varphi\in\mathcal{M}$,
and $\sigma<s$, the following a~priori estimate holds true:
$$
\|u\|_{H^{s,\varphi}(\Gamma)}\leq
c\,\bigr(\,\|Au\|_{H^{s-r,\,\varphi}(\Gamma)}+
\|u\|_{H^{\sigma,\varphi}(\Gamma)}\,\bigl)\quad\forall\;\;u\in
H^{s,\varphi}(\Gamma).
$$
Here the number $c>0$ does not depend on $u$.
\end{theorem}

If the spaces $N$ and $N^{+}$ are trivial, then the operator (2.1) is a
topological isomorphism. Generally, it is convenient to construct the
isomorphism with the help of two projectors. Let us decompose the spaces from
(2.1) into the following direct sums of (closed) subspaces:
$$
H^{s,\varphi}(\Gamma)=N\dotplus\bigl\{u\in H^{s,\varphi}(\Gamma):\,
(u,w)_{\Gamma}=0\;\;\forall\;\;w\in N\bigr\},
$$
$$
H^{s-r,\,\varphi}(\Gamma)=N^{+}\dotplus\bigl\{f\in
H^{s-r,\,\varphi}(\Gamma):\,(f,v)_{\Gamma}=0\;\;\forall\;\; v\in N^{+}\bigr\}.
$$
We denote by $P$ and $P^{+}$ respectively the projectors of these spaces on the
second terms in the sums in parallel to the first terms. The projectors do not
depend on $s$, $\varphi$.

\begin{theorem} 
Let $s\in\mathbb{R}$ and $\varphi\in\mathcal{M}$. The restriction of the
operator $(2.1)$ to the subspace $P(H^{s,\varphi}(\Gamma))$ establishes the
topological isomorphism
$$
A:\,P(H^{s,\varphi}(\Gamma))\leftrightarrow P^{+}(H^{s-r,\,\varphi}(\Gamma)).
$$
\end{theorem}

Theorems 2.1--2.3 were proved in \cite[Sec. 4]{M10}. They specify, with regard
to the refined scale, the known theorems on properties of an elliptic
pseudodifferential operator on the Sobolev scale (see \cite[Theorem
19.2.1]{Hermander85} or \cite[Theorems 2.3.3, 2.3.12]{Agranovich94}). Note that
the boundedness of the operator (2.1) holds true without the assumption about
ellipticity of $A$. If $\dim\Gamma\geq2$, then the index of operator (2.1) is
equal to zero \cite{AtSin63}, \cite[Sec. 2.3 f]{Agranovich94}. In the case
where $\dim\Gamma=1$, the index can be nonzero. There is a class of elliptic
operators depending on a complex parameter (so called parameter elliptic
operators) such that $N=N^{+}=\{0\}$ for all values of the parameter
sufficiently large in modulus \cite[Sec. 4.1]{Agranovich94}. Moreover for a
solution to a parameter elliptic equation, a certain two-sided a priory
estimate holds with constants independent of the parameter. Such an estimate
was obtained for the refined scale in \cite[Theorem 6.1]{M10}. The analogs of
Theorems 2.1--2.3 for different types of elliptic matrix operators were proved
in \cite{M8, M13, MM14}.

Let us study a local smoothness of an elliptic equation solution in the refined
scale. Let $\Gamma_{0}$ be an nonempty open set on the manifold $\Gamma$. We
denote
$$
H^{s,\varphi}_{\mathrm{loc}}(\Gamma_{0})
:=\bigl\{f\in\mathcal{D}'(\Gamma):\,\chi\,f\in
H^{s,\varphi}(\Gamma)\;\;\forall\;\;\chi\in
C^{\infty}(\Gamma),\;\mathrm{supp}\,\chi\subseteq \Gamma_{0}\bigr\}.
$$

\begin{theorem} 
Let $u\in\mathcal{D}'(\Gamma)$ be a solution to the equation $Au=f$ on
$\Gamma_{0}$ with $f\in H^{s,\varphi}_{\mathrm{loc}}(\Gamma_{0})$ for some
$s\in\mathbb{R}$ and $\varphi\in\mathcal{M}$. Then $u\in
H^{s+r,\,\varphi}_{\mathrm{loc}}(\Gamma_{0})$.
\end{theorem}

This theorem and the analog of Theorem 1.1 (vi) for the refined scale over
$\Gamma$ imply the following sufficient condition for a solution $u$ to have
continuous derivatives of a prescribed order.

\begin{theorem} 
Let $u\in\mathcal{D}'(\Gamma)$ be a solution to the equation $Au=f$ on
$\Gamma_{0}$, where $f\in H^{k-r+n/2,\,\varphi}_{\mathrm{loc}}(\Gamma_{0})$ for
a certain integer $k\geq0$ and a function parameter $\varphi$ satisfying
inequality $(1.5)$. Then $u\in C^{k}(\Gamma_{0})$.
\end{theorem}

Theorems 2.4 and 2.5 were proved in \cite[Sec. 5]{M10}. Theorem 2.5 shows an
advantage of the refined scale over the Sobolev scale when a classical
smoothness of a solution is under investigation. Indeed, if we restrict
ourselves to the case of $\varphi\equiv1$, we have to replace the condition
$f\in H^{k-r+n/2,\,\varphi}_{\mathrm{loc}}(\Gamma_{0})$ with the condition
$f\in H^{k-r+\varepsilon+n/2,\,1}_{\mathrm{loc}}(\Gamma_{0})$ for some
$\varepsilon>0$. The last condition is far stronger than previous one. The
analogs of Theorems 2.4 and 2.5 for elliptic matrix operators were proved in
\cite{M8, M13, MM14}. A local regularity of an elliptic system solution in the
Sobolev scale was investigated in \cite[Sec. 10.6]{Hermander63}. We also note
that, in the H\"ormander spaces, regularity properties of solutions to
hypoelliptic partial differential equations with constant coefficients were
studied  in \cite[Ch. IV]{Hermander63}, \cite[Ch.~11]{Hermander83}

At the end of this section we give, with the help of $A$, an alternative and
equivalent definition of the refined scale over the closed manifold $\Gamma$.

Let us assume that $\mathrm{ord}\,A=r>0$ and that the operator
$A:C^{\infty}(\Gamma)\rightarrow C^{\infty}(\Gamma)$ is positive in the space
$L_{2}(\Gamma,dx)$. We denote by $A_{0}$ the closure of this operator in
$L_{2}(\Gamma,dx)$. Let $s\in\mathbb{R}$, $\varphi\in\mathcal{M}$, and
$$
\varphi_{s,r}(t):=t^{s/r}\varphi(t^{1/r})\;\;\mbox{for}\;\;t\geq1\quad\mbox{and}
\quad\varphi_{s,r}(t):=\varphi(1)\;\;\mbox{for}\;\;0<t<1.
$$
The operator $\varphi_{s,r}(A_{0})$ is regarded in $L_{2}(\Gamma,dx)$ as the
Borel function $\varphi_{s,r}$ of the self-adjoint operator~$A_{0}$.

\begin{theorem} 
For arbitrary $s\in\mathbb{R}$ and $\varphi\in\mathcal{M}$, the space
$H^{s,\varphi}(\Gamma)$ coincides with the completion of the set of all
functions $u\in C^{\infty}(\Gamma)$ with respect to the norm
$\|\varphi_{s,r}(A_{0})\,u\|_{L_{2}(\Gamma)}$, which is equivalent to the norm
$\|u\|_{H^{s,\varphi}(\Gamma)}$.
\end{theorem}

An important example of the operator $A$ mentioned above is the operator
$1-\triangle_{\Gamma}$, where $\triangle_{\Gamma}$ is the Beltrami-Laplace
operator on the Riemannian manifold $\Gamma$ (then $r=2$).

Theorem 2.6 was proved in \cite[Sec 3.8]{MM11}. For equivalent definition of
the Sobolev scale over $\Gamma$, the powers of $A_{0}$ is used  instead of the
regular varying function $\varphi_{s}$ (see \cite[Sec 5.3]{Agranovich94}).

\section{An elliptic boundary problem on the one-sided scale} 

Let us recall that $\Omega$ is a bounded domain in $\mathbb{R}^n$, were
$n\geq2$, and that its boundary $\partial\Omega$ is a closed infinitely smooth
manifold of the dimension $n-1$. We consider the nonhomogeneous boundary
problem in the domain $\Omega$:
$$
L\,u\equiv\sum_{|\mu|\leq2q}\,l_{\mu}\,D^{\mu}u=f\;\;\mbox{in}\;\;\Omega,
\eqno(3.1)
$$
$$
B_{j}\,u\equiv\sum_{|\mu|\leq m_{j}}\,b_{j,\mu}\,D^{\mu}u=
g_{j}\;\;\mbox{on}\;\;\partial\Omega,\;\;j=1,\ldots,q. \eqno(3.2)
$$
Here $L$ and $B_{j}$ are linear partial differential expressions with
complex-valued coefficients $l_{\mu}\in C^{\infty}(\,\overline{\Omega}\,)$ and
$b_{j,\mu}\in C^{\infty}(\partial\Omega)$. We suppose that $\mathrm{ord}\,L=2q$
is an even positive number and $\mathrm{ord}\,B_{j}=m_{j}\leq2q-1$ for all
$j=1,\ldots,q$. Let $m:=\max\,\{m_{1},\ldots,m_{q}\}$.

In what follows the boundary problem (3.1), (3.2) is assumed to be
\textit{regular elliptic}. It means that the expression $L$ is proper elliptic
in $\overline{\Omega}$, and the system $B:=(B_{1},\ldots,B_{q})$ of the
boundary expressions is normal and satisfies the complementing condition with
respect to $L$ on $\partial\Omega$ (see \cite{LM72}, \cite[Sec.
5.2.1]{Triebel78}). It follows from the condition of normality that all numbers
$m_{j}$, $j=1,\ldots,q$, are distinct.

We will investigate the mapping $u\mapsto(Lu,Bu)$ in appropriate spaces of the
refined scales. To describe the range of this mapping, we consider the boundary
problem
$$
L^{+}v=\omega\;\;\mbox{in}\;\;\Omega, \eqno(3.3)
$$
$$
B^{+}_{j}v=h_{j}\;\;\mbox{on}\;\;\partial\Omega,\;\;j=1,\ldots,q, \eqno(3.4)
$$
formally adjoint to the problem (3.1), (3.2) with respect to the Green formula
$$
(Lu,v)_{\Omega}+\sum_{j=1}^{q}\;(B_{j}u,\,C_{j}^{+}v)_{\partial\Omega}
=(u,L^{+}v)_{\Omega}+\sum_{j=1}^{q}\;(C_{j}u,\,B_{j}^{+}v)_{\partial\Omega},
\;\;u,v\in C^{\infty}(\,\overline{\Omega}\,). \eqno(3.5)
$$
Here $L^{+}$ is the linear differential expression formally adjoint to $L$, and
$\{B^{+}_{j}\}$, $\{C_{j}\}$, $\{C^{+}_{j}\}$ are some normal systems of linear
differential boundary expressions. Their coefficients are infinitely smooth,
and their orders satisfy the equalities
$$
\mathrm{ord}\,L^{+}=2q,\quad\mathrm{ord}\,B_{j}+\mathrm{ord}\,C^{+}_{j}=
\mathrm{ord}\,C_{j}+\mathrm{ord}\,B^{+}_{j}=2q-1.
$$
We denote $m_{j}^{+}:=\mathrm{ord}\,B_{j}^{+}$. In (3.5) and bellow, the
notations $(\cdot,\cdot)_{\Omega}$ and $(\cdot,\cdot)_{\partial\Omega}$ stand
for the inner products in the spaces $L_{2}(\Omega)$ and
$L_{2}(\partial\Omega)$ respectively, and also denote the extensions by
continuity of these products.

We set
$$
\mathcal{N}:=\{u\in
C^{\infty}(\,\overline{\Omega}\,):\;Lu=0\;\;\mbox{in}\;\;\Omega,\;\;
B_{j}u=0\;\;\mbox{on}\;\;\partial\Omega\;\;\forall\;\;j=1,\ldots,q\},
$$
$$
\mathcal{N}^{+}:=\{v\in
C^{\infty}(\,\overline{\Omega}\,):\;L^{+}v=0\;\;\mbox{in}\;\; \Omega,\;\;
B^{+}_{j}v=0\;\;\mbox{on}\;\;\partial\Omega\;\;\forall\;\;j=1,\ldots,q\}.
$$
Since both the problems (3.1), (3.2) and (3.3), (3.4) are regular elliptic,
both the spaces $\mathcal{N}$ and $\mathcal{N}^{+}$ are finite dimensional.

\begin{theorem} 
Let $s>m+1/2$ and $\varphi\in\mathcal{M}$. The mapping
$$
(L,B):\,u\rightarrow(Lu,B_{1}u,\ldots,B_{q}u),\quad u\in
C^{\infty}(\,\overline{\Omega}\,), \eqno(3.6)
$$
is extended by a continuity to the bounded linear operator
$$
(L,B):\,H^{s,\varphi}(\Omega)\rightarrow
H^{s-2q,\,\varphi}(\Omega)\oplus\bigoplus_{j=1}^{q}\,
H^{s-m_{j}-1/2,\,\varphi}(\partial\Omega)=:
\mathcal{H}_{s,\varphi}(\Omega,\partial\Omega). \eqno(3.7)
$$
This operator is a Fredholm one. Its kernel coincides with $\mathcal{N}$, and
its range is equal to the set
$$
\Bigl\{(f,g_{1},\ldots,g_{q})\in\mathcal{H}_{s,\varphi}(\Omega,\partial\Omega):\,
(f,v)_{\Omega}+\sum_{j=1}^{q}\,(g_{j},C^{+}_{j}v)_{\partial\Omega}=0\;\;
\forall\;\;v\in \mathcal{N}^{+}\Bigr\}. \eqno(3.8)
$$
The index of the operator \rm (3.7) \it is equal to
$\dim\mathcal{N}-\dim\mathcal{N}^{+}$ and does not depend on $s$, $\varphi$.
\end{theorem}

In this theorem and in the next theorems of the section, the condition
$s>m+1/2$ is essential. Indeed, if $s<m_{j}+1/2$ for some $j=1,\ldots,q$, then
the mapping $u\rightarrow B_{j}u$, $u\in C^{\infty}(\,\overline{\Omega}\,)$,
can not be extended to the continuous linear operator
$B_{j}:H^{s,\varphi}(\Omega)\rightarrow\mathcal{D}'(\partial\Omega)$. Thus the
operator (3.6) is correctly defined on the upper refined one-sided scale
$$
\{H^{s,\varphi}(\Omega):s>m+1/2,\varphi\in\mathcal{M}\,\bigr\}.
$$
Hence the left-hand sides of equations (3.1), (3.2) is defined for each $u\in
H^{s,\varphi}(\Omega)$ with $s>m+1/2$, whereas these equations are understood
in the theory of distributions.

\begin{theorem} 
For arbitrarily chosen parameters $s>m+1/2$, $\varphi\in\mathcal{M}$, and
$\sigma<s$, the following a~priori estimate holds true:
$$
\|u\|_{H^{s,\varphi}(\Omega)}\leq
c\,\bigr(\,\|(L,B)u\|_{\mathcal{H}_{s,\varphi}(\Omega,\partial\Omega)}+
\|u\|_{H^{\sigma,\varphi}(\Omega)}\,\bigl)\quad\forall\;\;u\in
H^{s,\varphi}(\Omega).
$$
Here the number $c>0$ does not depend on $u$.
\end{theorem}

If the spaces $\mathcal{N}$ and $\mathcal{N}^{+}$ are trivial, then the
operator (3.7) is a topological isomorphism. In general, we can get the
isomorphism with the help of two projectors. Let the spaces in which the
operator (3.7) acts be decomposed into the following direct sums of subspaces:
$$
H^{s,\varphi}(\Omega)=\mathcal{N}\dotplus\bigl\{u\in
H^{s,\varphi}(\Omega):\;(u,w)_{\Omega}=0\;\;\forall\;\;w\in\mathcal{N} \bigr\},
$$
$$
\mathcal{H}_{s,\varphi}(\Omega,\partial\Omega)=\bigl\{(v,0,\ldots,0):v\in
\mathcal{N}^{+}\bigr\}\dotplus\,(3.8).
$$
We denote by $\mathcal{P}$ and $\mathcal{Q}^{+}$ respectively the projectors of
these spaces on the second terms in the sums in parallel to the first terms.
The projectors are independent of $s$ and $\varphi$.

\begin{theorem} 
Let $s>m+1/2$ and $\varphi\in\mathcal{M}$. The restriction of the operator
$(3.7)$ to the subspace $\mathcal{P}(H^{s,\varphi}(\Omega))$ establishes the
topological isomorphism
$$
(L,B):\,\mathcal{P}(H^{s,\varphi}(\Omega))\leftrightarrow
\mathcal{Q}^{+}(\mathcal{H}_{s,\varphi}(\Omega,\partial\Omega)).
$$
\end{theorem}

Theorems 3.1--3.3 were proved in \cite[Sec. 4]{MM2b}. The boundedness of the
operator (3.7) holds true without the assumption that the boundary problem
(3.1), (3.2) is elliptic. In the paper \cite{Slenzak74} this problem was
studied in a different scale of the H\"ormander spaces (also called a refined
one). Theorems 3.1--3.3 specify, with regard to the refined scale, the known
theorems on properties of an elliptic boundary problem in the Sobolev one-sided
scale (see \cite[Ch. V]{ADN59}, \cite[Ch. 2, Sec. 5.4]{LM72}, \cite[Ch.
20]{Hermander85}, \cite[Sec 2, 4]{Agranovich97}).  The analogs of Theorems
3.1--3.3 are valid for nonregular elliptic boundary problems \cite{MM2b} and
for elliptic problems for systems of partial differential equations \cite{M9}.
The case where the boundary operators have distinct orders on different
connected components of the domain $\Omega$ was considered especially in
\cite{M7}. There is a class of elliptic boundary problems depending on a
parameter $\lambda\in\mathbb{C}$ such that $\mathcal{N}=\mathcal{N}^{+}=\{0\}$
for $|\lambda|\gg1$, and hence the index of the corresponding operator is equal
to $0$ for all $\lambda$ (see \cite{AgmNir63, AgrVish64},
\cite[Sec.~3]{Agranovich97}). For a solution to such a parameter elliptic
problem, a certain two-sided a priory estimate holds with constants independent
of the parameter $\lambda\in\mathbb{C}$ with $|\lambda|\gg1$. Such an estimate
was obtained for the refined scale in \cite[Theorem 7.2]{MM2c}. Regular
elliptic boundary problems in positive one-sided scales of different normed
spaces were studied in \cite{ADN59, Triebel78, Triebel83}.

Now we study an increase in a local smoothness of an elliptic boundary problem
solution. Let $U$ be an open subset in $\mathbb{R}^{n}$. We set
$\Omega_{0}:=U\cap\Omega\neq\emptyset$ and $\Gamma_{0}:=U\cap\partial\Omega$
(the case were $\Gamma_{0}=\emptyset$ is possible). Let us introduce the
following local analogs of spaces of the refined scales:
$$
H^{\sigma,\varphi}_{\mathrm{loc}}(\Omega_{0},\Gamma_{0}):=
\bigl\{u\in\mathcal{D}'(\Omega):\chi\,u\in
H^{\sigma,\varphi}(\Omega)\;\;\forall\;\;\chi\in
C^{\infty}(\overline{\Omega}),\;\mathrm{supp}\,\chi\subseteq\Omega_{0}\cup\Gamma_{0}
\bigr\},
$$
$$
H^{\sigma,\varphi}_{\mathrm{loc}}(\Gamma_{0}):=\bigl\{h\in
\mathcal{D}'(\partial\Omega): \chi\,h\in
H^{\sigma,\varphi}(\partial\Omega)\;\;\forall\;\;\chi\in
C^{\infty}(\partial\Omega),\;\mathrm{supp}\,\chi\subseteq\Gamma_{0}\bigr\}.
$$
Here $\sigma\in\mathbb{R}$, $\varphi\in\mathcal{M}$ and, as usual,
$\mathcal{D}'(\Omega)$ denotes the topological space of all distributions in
$\Omega$.

\begin{theorem} 
Let  $s>m+1/2$ and $\eta\in\mathcal{M}$. Suppose that the distribution $u\in
H^{s,\eta}(\Omega)$ is a solution to the problem $(3.1)$, $(3.2)$, where
$$
f\in H^{s-2q+\varepsilon,\,\varphi}_{\mathrm{loc}}(\Omega_{0},\Gamma_{0})
\quad\mbox{and}\quad g_{j}\in
H^{s-m_{j}-1/2+\varepsilon,\,\varphi}_{\mathrm{loc}}(\Gamma_{0}),
\;\;j=1,\ldots,q,
$$
for some $\varepsilon\geq0$ and $\varphi\in\mathcal{M}$. Then $u\in
H^{s+\varepsilon,\,\varphi}_{\mathrm{loc}}(\Omega_{0},\Gamma_{0})$.
\end{theorem}

Note that in the case where $\Omega_{0}=\Omega$ and $\Gamma_{0}=\partial\Omega$
we have the global smoothness increase (i.e. the increase in the whole closed
domain $\overline{\Omega}$). If $\Gamma_{0}=\emptyset$, then we get an interior
smoothness increase (in an open subset $\Omega_{0}\subseteq\Omega$).

Theorems 3.4 and 1.1 (vi) imply the following sufficient condition for the
solution $u$ to be classical.

\begin{theorem} 
Let  $s>m+1/2$ and $\chi\in\mathcal{M}$. Suppose that the distribution $u\in
H^{s,\chi}(\Omega)$ is a solution to the problem $(3.1)$, $(3.2)$ in which
$$
f\in H^{n/2,\,\varphi}_{\mathrm{loc}}(\Omega,\emptyset)\cap
H^{m-2q+n/2,\,\varphi}(\Omega),
$$
$$
g_{j}\in H^{m-m_{j}+(n-1)/2,\,\varphi}(\partial\Omega),\;\;j=1,\ldots,q,
$$
and the function parameter $\varphi\in\mathcal{M}$ satisfies condition $(1.5)$.
Then the solution $u$ is classical, that is $u\in C^{2q}(\Omega)\cap
C^{m}(\,\overline{\Omega}\,)$.
\end{theorem}

Theorems 3.4, 3.5 were proved in \cite[Sec. 5, 6]{MM2c} (generally, for a non
regular elliptic problem). The analog of Theorem 3.4 is valid for elliptic
boundary problems for systems of partial differential equations \cite{M9}. In
the Sobolev positive one-sided scale ($s\geq0$, $\varphi\equiv1$), a smoothness
of solutions to elliptic boundary problems was investigated in
\cite{Nirenberg55, Browder56, Schechter61}, \cite[Ch. 3, Sec. 4]{Berezansky68}
(see also \cite[Sec. 2.4]{Agranovich97}).

\section{Semihomogeneous elliptic problems} 

\subsection{} 
As we have mentioned, the results of Section 3 are not valid for $s<m+1/2$
because the mapping (3.6) can not be extended to the bounded linear operator
(3.7). But if the boundary problem (3.1), (3.2) is semihomogeneous (i.e.,
$f\equiv0$ or all $g_{j}\equiv0$), it establishes a bounded and Fredholm
operator in the two-sided refined scale (for all real $s$). We will consider
separately the case of the homogeneous elliptic equation (3.1) and the case of
the homogeneous boundary conditions (3.2).

\subsection{A boundary problem for a homogeneous elliptic equation} 
Let us consider the regular elliptic boundary problem (3.1), (3.2), provided
that $f\equiv0$:
$$
Lu=0\;\;\mbox{on}\;\;\Omega,\quad
B_{j}u=g_{j}\;\;\mbox{on}\;\;\partial\Omega,\;\;j=1,\ldots,q. \eqno(4.1)
$$

We will connect the following spaces with this problem:
$$
K_{L}^{\infty}(\Omega):=\bigl\{\,u\in
C^{\infty}(\,\overline{\Omega}\,):\,L\,u=0\;\;\mbox{in}\;\;\Omega\,\bigr\},
$$
$$
K_{L}^{s,\varphi}(\Omega):=\bigl\{\,u\in
H^{s,\varphi}(\Omega):\,L\,u=0\;\;\mbox{in}\;\;\Omega\,\bigr\}
$$
for $s\in\mathbb{R}$, $\varphi\in\mathcal{M}$. It  follows from a continuity of
the embedding $H^{s,\varphi}(\Omega)\hookrightarrow\mathcal{D}'(\Omega)$ that
$K_{L}^{s,\varphi}(\Omega)$ is a closed subspace in $H^{s,\varphi}(\Omega)$. We
can consider $K_{L}^{s,\varphi}(\Omega)$ as a Hilbert space with respect to the
inner product in $H^{s,\varphi}(\Omega)$.

\begin{theorem} 
Let $s\in\mathbb{R}$ and $\varphi\in\mathcal{M}$. The set
$K_{L}^{\infty}(\Omega)$ is dense in the space $K_{L}^{s,\varphi}(\Omega)$. The
mapping
$$
u\mapsto Bu=(B_{1}u,\ldots,B_{q}u),\;\;u\in K_{L}^{\infty}(\Omega),
$$
is extended by a continuity to the bounded linear operator
$$
B:\,K_{L}^{s,\varphi}(\Omega)\rightarrow
\bigoplus_{j=1}^{q}\,H^{s-m_{j}-1/2,\,\varphi}(\partial\Omega)=:
\mathcal{H}_{s,\varphi}(\partial\Omega). \eqno(4.2)
$$
This operator is a Fredholm one. Its kernel coincides with $\mathcal{N}$, and
its range is equal to the set
$$
\Bigl\{(g_{1},\ldots,g_{q})\in\mathcal{H}_{s,\varphi}(\partial\Omega):\,
\sum_{j=1}^{q}\,(g_{j},C^{+}_{j}v)_{\partial\Omega}=0\;\;\forall\;\;v\in
\mathcal{N}^{+}\Bigr\}.
$$
The index of the operator \rm (4.2) \it is equal to
$\dim\mathcal{N}-\dim\mathcal{G}^{+}$ where
$$
\mathcal{G}^{+}:=\bigl\{\,\bigl(C_{1}^{+}v,\ldots,C_{q}^{+}v\bigr):\,v\in
\mathcal{N}^{+}\,\bigr\},
$$
and does not depend on $s$, $\varphi$.
\end{theorem}

Theorem 4.1 was proved in \cite[Sec. 6]{MM5}. In contrast to Theorem 3.1, the
ellipticity condition is essential for the boundedness of the operator (4.2) in
the case where $s\leq m+1/2$. Note that $\dim\mathcal{G}^{+}
\leq\dim\mathcal{N}^{+}$ where the strict inequality is possible that results
from \cite[Theorem 13.6.15]{Hermander83}. In the case where $\varphi\equiv1$
and $s\in\mathbb{R}\setminus\{-1/2,-3/2,-5/2,\ldots\}$ Theorem 4.1 is a
consequence of the Lions--Magenes Theorems \cite[Ch. 2, Sec. 6.6, 7.3]{LM72}
(see also \cite{LM62V, LM63VI} and \cite[Sec. 6.10, 6.12]{Magenes65}).

\subsection{An elliptic problem with homogeneous boundary conditions} 
Now we will consider the regular elliptic boundary problem (3.1), (3.2),
provided that all $g_{j}\equiv0$:
$$
Lu=f\;\;\mbox{in}\;\;\Omega,\quad
B_{j}u=0\;\;\mbox{on}\;\;\partial\Omega,\;\;j=1,\ldots,q. \eqno(4.3)
$$

Let us introduce the function spaces in which the operator of the problem (4.3)
acts. For the sake of brevity, we denote by $(\mathrm{b.c.})$ the homogeneous
boundary conditions in (4.3). In addition, we denote by $(\mathrm{b.c.})^{+}$
the homogeneous boundary conditions
$$
B^{+}_{j}v=0\;\;\mbox{on}\;\;\partial\Omega,\;\;j=1,\ldots,q.
$$
They correspond to the formally adjoint boundary problem (3.3), (3.4). We set
$$
C^{\infty}(\mathrm{b.c.}):=\bigl\{u\in
C^{\infty}(\,\overline{\Omega}\,):\,B_{j}u=0\;\;\mbox{on}\;\;\partial\Omega\;\;
\forall\;\;j=1,\ldots,q\bigr\},
$$
$$
C^{\infty}(\mathrm{b.c.})^{+}:=\bigl\{v\in
C^{\infty}(\,\overline{\Omega}\,):\,B^{+}_{j}v=0\;\;\mbox{on}\;\;\partial\Omega
\;\;\forall\;\;j=1,\ldots,q\bigr\}.
$$

Let $s\in\mathbb{R}$ and $\varphi\in\mathcal{M}$. We define the Hilbert space
$H^{s,\varphi,(0)}(\Omega)$ in the following way:
$$
H^{s,\varphi,(0)}(\Omega):=
\begin{cases}
\;H^{s,\varphi}(\Omega)\;\; & \text{for}\;\;s\geq0, \\
\;H^{s,\varphi}_{\overline{\Omega}}(\mathbb{R}^{n}) & \text{for}\;\;s<0.
\end{cases}
$$
According to Theorem 1.1 (iii), (vii), the spaces $H^{s,\varphi,(0)}(\Omega)$
and $H^{-s,1/\varphi,(0)}(\Omega)$ are mutually dual for every $s\in\mathbb{R}$
with respect to the inner product in $L_{2}(\Omega)$. It also follows from
Theorem 1.1 (i), (ii) that the set $C^{\infty}(\,\overline{\Omega}\,)$ is dense
in the space $H^{s,\varphi,(0)}(\Omega)$ for each $s\in\mathbb{R}$. Here we
identify each function $f\in C^{\infty}(\,\overline{\Omega}\,)$ with its
extension by zero
$$
\mathcal{O}f(x):=
\begin{cases}
\;f(x) &\;\; \text{for}\;\; x\in\overline{\Omega}, \\
\;0 &\;\; \text{for}\;\; x\in\mathbb{R}^{n}\setminus\overline{\Omega},
\end{cases} \eqno(4.4)
$$
which is a regular distribution in
$H^{s,\varphi}_{\overline{\Omega}}(\mathbb{R}^{n})$ for $s<0$. Now one may
conclude that Theorem 1.1 (iii), (iv) implies the continuous dense embedding
$$
H^{s_{1},\varphi_{1},(0)}(\Omega)\hookrightarrow
H^{s,\varphi,(0)}(\Omega)\quad\mbox{for}\;\;-\infty<s<s_{1}<\infty,\;\;\mbox{and}
\;\;\varphi,\varphi_{1}\in\mathcal{M}.
$$

Finally, let us define the Hilbert spaces $H^{s,\varphi}(\mathrm{b.c.})$ and
$H^{s,\varphi}(\mathrm{b.c.})^{+}$ of distributions satisfying the homogeneous
boundary conditions. In the case where $s\notin\{m_{j}+1/2:j=1,\ldots,q\}$ we
denote by $H^{s,\varphi}(\mathrm{b.c.})$ the closure of
$C^{\infty}(\mathrm{b.c.})$ in the space $H^{s,\varphi,(0)}(\Omega)$. In the
case where $s\in\{m_{j}+1/2:j=1,\ldots,q\}$ we define the space
$H^{s,\varphi}(\mathrm{b.c.})$ by means of the interpolation with the parameter
$\psi(t)=t^{1/2}$:
$$
H^{s,\varphi}(\mathrm{b.c.}):=\bigl[H^{s-1/2,\,\varphi}(\mathrm{b.c.}),
H^{s+1/2,\,\varphi}(\mathrm{b.c.})\bigr]_{t^{1/2}}. \eqno(4.5)
$$
If we change $(\mathrm{b.c.})$ for $(\mathrm{b.c.})^{+}$, and $m_{j}$ for
$m_{j}^{+}$ in the last two sentences, we give the definition of the space
$H^{s,\varphi}(\mathrm{b.c.})^{+}$. Note that in the case where
$s\in\{m_{j}+1/2:j=1,\ldots,q\}$ the norms in the spaces
$H^{s,\varphi}(\mathrm{b.c.})$ and $H^{s,\varphi,(0)}(\Omega)$ are not
equivalent. The analogous fact is true for $H^{s,\varphi}(\mathrm{b.c.})^{+}$.

\begin{proposition} 
Let $s>0$, $s\neq m_{j}+1/2$ for all $j=1,\ldots,q$, and
$\varphi\in\mathcal{M}$. Then
$$
H^{s,\varphi}(\mathrm{b.c.})=\bigl\{u\in
H^{s,\varphi}(\Omega):B_{j}u=0\;\mbox{on}\;\partial\Omega\;\mbox{for
all}\;j=1,\ldots,q\;\mbox{such that}\;\;s>m_{j}+1/2\bigr\}.
$$
If $s<1/2$, then $H^{s,\varphi}(\mathrm{b.c.})=H^{s,\varphi,(0)}(\Omega)$. This
proposition remains true if we change $m_{j}$ for $m_{j}^{+}$,
$(\mathrm{b.c.})$ for $(\mathrm{b.c.})^{+}$, and $B_{j}$ for $B_{j}^{+}$.
\end{proposition}

\begin{theorem} 
Let $s\in\mathbb{R}$ and $\varphi\in\mathcal{M}$. The mapping $u\mapsto Lu$,
$u\in C^{\infty}(\mathrm{b.c.})$, is extended by a continuity to the bounded
linear operator
$$
L:H^{s,\varphi}(\mathrm{b.c.})\rightarrow
(H^{2q-s,\,1/\varphi}(\mathrm{b.c.})^{+})'. \eqno(4.6)
$$
Here the function $Lu$ is interpreted as the functional
$(Lu,\,\cdot\,)_{\Omega}$, whereas $(H^{2q-s,\,1/\varphi}(\mathrm{b.c.})^{+})'$
denotes the antidual space to $H^{2q-s,\,1/\varphi}(\mathrm{b.c.})^{+}$ with
respect to the inner product in $L_{2}(\Omega)$. The operator $(4.6)$ is a
Fredholm one. Its kernel coincides with $\mathcal{N}$, and its range is equal
to the set
$$
\bigl\{\,f\in(H^{2q-s,\,1/\varphi}(\mathrm{b.c.})^{+})':\,
(f,v)_{\Omega}=0\;\;\forall\;\;v\in\mathcal{N}^{+}\,\bigr\}.
$$
The index of the operator $(4.6)$ is equal to
$\dim\mathcal{N}-\dim\mathcal{N}^{+}$ and does not depend on $s$, $\varphi$.
\end{theorem}

Theorem 4.2 was proved in \cite[Sec. 5]{MM6}, provided that $s\neq j-1/2$ for
each $j=1,\ldots,2q$. For the rest values of $s$, the theorem is deduced by
means of the interpolation formula (4.5). The analogs of Theorems 3.2--3.4 was
obtained for the operator (4.6) as well (see also \cite{MM1}). Theorem 4.2
specifies, with regard to the refined scale, the theorem of Berezansky, Krein
and Roitberg on homeomorphisms realized by the elliptic operator $L$ on the
two-sided Sobolev scale \cite{BKR63}, \cite[Ch. 3, Sec. 6]{Berezansky68},
\cite[Sec. 5.5]{Roitberg96}. In the case of $s\leq m+1/2$ the ellipticity
condition is essential for the boundedness of the operator (4.6). The
interpolation space (4.5) was studied in the Sobolev case of $\varphi\equiv1$
in \cite{Grisvard67, Seeley72} (see also \cite[Sec. 4.3.3]{Triebel78}).

\subsection{} 
We note that the general nonhomogeneous boundary problem (3.1), (3.2) cannot be
reduced to the semihomogeneous boundary problems in the lower part of the
refined scale, namely for $s<m+1/2$. Indeed, if $s<-1/2$, then solutions to
these problems belong to the spaces of distributions of the different nature;
solutions to the problem (4.1) belong to $K_{L}^{s,\varphi}(\Omega)\subset
H^{s,\varphi}(\Omega)$ being distributions defined in the open domain $\Omega$,
whereas solutions to the problem (4.3) belong to
$H^{s,\varphi}(\mathrm{b.c.})\subset
H^{s,\varphi}_{\overline{\Omega}}(\mathbb{R}^{n})$ being distributions
supported on the closed domain $\overline{\Omega}$. If $-1/2<s<m+1/2$, then
solutions to the semihomogeneous problems are distributions defined in $\Omega$
(see Theorem 1.1 (iii) in the case $-1/2<s<0$), but the operator $(L,B)$ can
not be correctly defined on $K_{L}^{s,\varphi}(\Omega)\cup
H^{s,\varphi}(\mathrm{b.c.})$ because of the inequality
$$
(K_{L}^{s,\varphi}(\Omega)\cap
H^{s,\varphi}(\mathrm{b.c.}))\setminus\mathcal{N}\neq\emptyset. \eqno(4.7)
$$

Note also that in the case where $s>m+1/2$ we have the equality of sets in
(4.7). Hence the nonhomogeneous problem (3.1), (3.2) is reduced to the
semihomogeneous problems (4.1) and (4.3); i.e., Theorem 3.1 is equivalent to
Theorems 4.1 and 4.2 taken together.

\section{Generic theorems for elliptic problems in two-sided scales}

In \cite{Roitberg64, Roitberg96, Roitberg99} Ya. A. Roitberg introduced a
special modification of the Sobolev two-sided scale in which the operator of an
elliptic boundary problem is bounded and a Fredholm one for every parameter
$s\in\mathbb{R}$  (see also \cite[Ch. 3, Sec. 6]{Berezansky68}, \cite[Sec.
7.9]{Agranovich97}). This modification does not depend on coefficients of the
elliptic differential expression but depends solely on the order of the
expression. Therefore, the theorems on properties of elliptic problems in such
modified scale is naturally to call generic (for the class of the problems
having the same order). We will consider these theorems with regard to the
refined scale.

Let $s\in\mathbb{R}$, $\varphi\in\mathcal{M}$, and integer $r>0$. We set
$E_{r}:=\{k-1/2:k=1,\ldots,r\}$. In the case where $s\in\mathbb{R}\setminus
E_{r}$ we denote by $H^{s,\varphi,(r)}(\Omega)$ the completion of
$C^{\infty}(\,\overline{\Omega}\,)$ with respect to the Hilbert norm
$$
\|u\|_{H^{s,\varphi,(r)}(\Omega)}:=
\Bigl(\,\|u\|_{H^{s,\varphi,(0)}(\Omega)}^{2}+
\sum_{k=1}^{r}\;\bigl\|(D_{\nu}^{k-1}u)\upharpoonright\partial\Omega\,\bigr\|
_{H^{s-k+1/2,\varphi}(\partial\Omega)}^{2}\,\Bigr)^{1/2}.
$$
Here $D_{\nu}:=i\,\partial/\partial\nu$, with $\nu$ being the unit vector of
the inner normal to $\partial\Omega$. In the case where $s\in E_{r}$ we set
$$
H^{s,\varphi,(r)}(\Omega):=\bigl[\,H^{s-1/2,\varphi,(r)}(\Omega),
H^{s+1/2,\varphi,(r)}(\Omega)\,\bigr]_{t^{1/2}}.
$$

The collection of separable Hilbert spaces
$$
\{H^{s,\varphi,(r)}(\Omega):s\in\mathbb{R},\varphi\in\mathcal{M}\,\} \eqno(5.1)
$$
is called the refined scale modified in the Roitberg sense. The number $r$ is
called the index of this modification.

The scale (5.1) admits the following description. Let us denote by
$\Upsilon_{s,\varphi,(r)}(\Omega,\partial\Omega)$ the space of all
vector-functions
$$
(u_{0},u_{1},\ldots,u_{r})\in H^{s,\varphi,(0)}(\Omega)
\oplus\bigoplus_{k=1}^{r}\,H^{s-k+1/2,\,\varphi}(\partial\Omega) \eqno(5.2)
$$
such that $u_{k}=(D_{\nu}^{k-1}u_{0})\upharpoonright\partial\Omega$ for every
integer $k=1,\ldots r$ satisfying the inequality $s>k-1/2$. In view of Theorem
1.1 (viii), $\Upsilon_{s,\varphi,(r)}(\Omega,\partial\Omega)$ is a Hilbert
space with respect to the inner product in the space (5.2).

\begin{proposition} 
The mapping
$$
T_{r}:u\mapsto\bigl(\,u,u\upharpoonright\partial\Omega,\ldots,
(D_{\nu}^{r-1}u)\upharpoonright\partial\Omega\,\bigr),\quad u\in
C^{\infty}(\,\overline{\Omega}\,),
$$
is extended by a continuity to the bounded linear injective operator
$$
T_{r}:\,H^{s,\varphi,(r)}(\Omega)\rightarrow
\Upsilon_{s,\varphi,(r)}(\Omega,\partial\Omega) \eqno(5.3)
$$
far all $s\in\mathbb{R}$ and $\varphi\in\mathcal{M}$. If $s\notin E_{r}$, then
the operator $(5.3)$ is an isometric isomorphism.
\end{proposition}

Thus, we can interpret an element $u\in H^{s,\varphi,(r)}(\Omega)$ as the
vector-valued function
$$
(u_{0},u_{1},\ldots,u_{r}):=T_{r}u\in
\Upsilon_{s,\varphi,(r)}(\Omega,\partial\Omega). \eqno(5.4)
$$
Note that in view of Theorem 1.1 (viii)
$$
\|u\|_{H^{s,\varphi,(r)}(\Omega)}\asymp\|u_{0}\|_{H^{s,\varphi,(0)}(\Omega)}=
\|u_{0}\|_{H^{s,\varphi}(\Omega)}\;\;\;\mbox{if}\;\;s>r-1/2.
$$
Therefore
$$
H^{s,\varphi,(r)}(\Omega)=H^{s,\varphi}(\Omega)\;\;\mbox{with equvivalence of
norms}\;\;\mbox{if}\;\;s>r-1/2. \eqno(5.5)
$$

\begin{theorem} 
Let $s\in\mathbb{R}$ and $\varphi\in\mathcal{M}$. The mapping $(3.6)$ is
extended by a continuity to the bounded linear operator
$$
(L,B):\,H^{s,\varphi,(2q)}(\Omega)\rightarrow
H^{s-2q,\varphi,(0)}(\Omega)\oplus\bigoplus_{j=1}^{q}\,
H^{s-m_{j}-1/2,\,\varphi}(\partial\Omega)=:
\mathcal{H}_{s,\varphi,(0)}(\Omega,\partial\Omega). \eqno(5.6)
$$
This operator is a Fredholm one. Its kernel coincides with $\mathcal{N}$, and
its range is equal to the set
$$
\Bigl\{(f,g_{1},\ldots,g_{q})\in\mathcal{H}_{s,\varphi,(0)}(\Omega,\partial\Omega):
\,(f,v)_{\Omega}+\sum_{j=1}^{q}\,(g_{j},C^{+}_{j}v)_{\partial\Omega}=0\;\;
\forall\;\;v\in \mathcal{N}^{+}\Bigr\}.
$$
The index of the operator \rm (5.6) \it is equal to
$\dim\mathcal{N}-\dim\mathcal{N}^{+}$ and does not depend on $s$, $\varphi$.
\end{theorem}

This theorem is generic because the spaces in which the operator (5.6) acts are
the same for all boundary problems of the common order
$(2q,m_{1},\ldots,m_{q})$. It follows from (5.5) that Theorem 5.1 coincides
with Theorem 3.1 for $s>2q-1/2$.

Using Proposition 5.1 we give the following interpretation of a solution $u\in
H^{s,\varphi,(2q)}(\Omega)$ to the boundary problem (3.1), (3.2) in the sense
of the distribution theory. Let us write down the differential expressions $L$
and $B_{j}$ in a neighborhood of $\partial\Omega$ in the form
$$
L=\sum_{k=0}^{2q}\;L_{k}\,D_{\nu}^{k},\quad
B_{j}=\sum_{k=0}^{m_{j}}\;B_{j,k}\,D_{\nu}^{k}. \eqno(5.7)
$$
Here $L_{k}$ and $B_{j,k}$ are certain tangent differential expression.
Integrating by parts we arrive at the (special) Green formula
$$
(Lu,v)_{\Omega}=(u,L^{+}v)_{\Omega}-
i\sum_{k=1}^{2q}\;(D_{\nu}^{k-1}u,L^{(k)}v)_{\partial\Omega},\quad u,v\in
C^{\infty}(\,\overline{\Omega}\,).
$$
Here $L^{(k)}:=\sum_{r=k}^{2q}D_{\nu}^{r-k}L_{r}^{+}$, with $L_{r}^{+}$ being
the tangent differential expression formally adjoint to $L_{r}$. By passing to
the limit and using the notation (5.4) we get the next equality for $u\in
H^{s,\varphi,(2q)}(\Omega)$:
$$
(Lu,v)_{\Omega}=(u_{0},L^{+}v)_{\Omega}-
i\sum_{k=1}^{2q}\;(u_{k},L^{(k)}v)_{\partial\Omega},\quad v\in
C^{\infty}(\,\overline{\Omega}\,). \eqno(5.8)
$$
Now it follows from (5.7), (5.8) that the element $u\in
H^{s,\varphi,(2q)}(\Omega)$ is a solution to the boundary problem (3.1), (3.2)
with $f\in H^{s-2q,\varphi,(0)}(\Omega)$, $g_{i}\in
H^{s-m_{j}-1/2,\,\varphi}(\partial\Omega)$ if and only if the following
equalities hold true:
$$
(u_{0},L^{+}v)_{\Omega}- i\sum_{k=1}^{2q}\;(u_{k},L^{(k)}v)_{\partial\Omega}=
(f,v)_{\Omega}\quad\mbox{for all}\quad v\in C^{\infty}(\,\overline{\Omega}\,),
$$
$$
\sum_{k=0}^{m_{j}}\;B_{j,k}\,u_{k+1}=g_{j}\;\;\mbox{on}\;\;\partial\Omega,\quad
j=1,\ldots,q.
$$

Theorem 5.1 was proved in \cite[Sec. 5]{MM12}. The analogs of Theorems 3.2--3.4
were obtained for the operator (5.6) as well. Theorem 5.1 specifies, with
regard to the refined scale, the theorem of Ya. A. Roitberg on the Fredholm
property of a regular elliptic boundary problem in the modified Sobolev scale
(so-called theorem on a complete collections of homeomorphisms)
\cite{Roitberg64}, \cite[Sec. 4.1, 5.3]{Roitberg96} (see also \cite[Ch. 3, Sec.
6]{Berezansky68}, \cite[Sec. 7.9]{Agranovich97}). The analogs of Theorem 5.1
are also valid for nonregular elliptic boundary problems both for one and for
system of partial differential equations. Note that the boundedness of the
operator (5.6) holds true without the ellipticity assumption. Elliptic boundary
problems in the modified two-sided scales of different normed spaces were
studied in \cite{Roitberg96} (the Sobolev $L_{p}$-spaces) and in
\cite{Murach94a, Murach94b} (non-Sobolev spaces). A certain classes of
non-elliptic problems were investigated in the two-sided modified scales as
well (see \cite{Roitberg99}, \cite{EidZhit98} and the references therein).

\section{Individual theorems for elliptic problems} 

In the individual theorems, the domain of the operator $(L,B)$ depends on
coefficients of the elliptic expression $L$. Namely, we consider the operator
$$
(L,B):\,D^{s,\varphi}_{L,X}(\Omega)\rightarrow
X(\Omega)\oplus\bigoplus_{j=1}^{q}\,
H^{s-m_{j}-1/2,\,\varphi}(\partial\Omega)=:
\mathcal{X}_{s,\varphi}(\Omega,\partial\Omega). \eqno(6.1)
$$
Here $s\in\mathbb{R}$, $\varphi\in\mathcal{M}$, and $X(\Omega)$ is a certain
Hilbert space consisting of distributions in $\Omega$ and satisfying the
continuous embedding $X(\Omega)\hookrightarrow\mathcal{D}'(\Omega)$. The domain
of the operator (6.1) is the Hilbert space
$$
D^{s,\varphi}_{L,X}(\Omega):=\bigl\{u\in H^{s,\varphi}(\Omega):\, Lu\in
X(\Omega)\bigr\}
$$
endowed with the graphics inner product
$$
(u,v)_{D^{s,\varphi}_{L,X}(\Omega)}:=
(u,v)_{H^{s,\varphi}(\Omega)}+(Lu,Lv)_{X(\Omega)}.
$$

In the case where $s>m+1/2$ we may set $X(\Omega):=H^{s-2q,\,\varphi}(\Omega)$
that leads us to Theorem~3.1. But in the case where $s\leq m+1/2$ we cannot do
so if we want to define the operator $(L,B)$ on the non-modified refined scale.
The space $X(\Omega)$ must be narrower than $H^{s-2q,\,\varphi}(\Omega)$.

Let us formulate the conditions on $X(\Omega)$ under which the operator (5.1)
is bounded and has the Fredholm property for some $s$ and $\varphi$.

\begin{condition} 
The set $X^{\infty}(\Omega):=X(\Omega)\cap C^{\infty}(\,\overline{\Omega}\,)$
is dense in the space $X(\Omega)$.
\end{condition}

\begin{condition} 
There exists a number $c>0$ such that
$$
\|\mathcal{O}f\|_{H^{s-2q,\varphi}(\mathbb{R}^{n})}\leq
c\,\|f\|_{X(\Omega)}\quad\forall\;\;f\in X^{\infty}(\Omega).
$$
\end{condition}

We recall that the function $\mathcal{O}f$ is given by formula (4.4). It
follows from the Conditions~1 and 2 that the mapping $f\mapsto\mathcal{O}f$,
$f\in X^{\infty}(\Omega)$, is extended by continuity to the linear bounded
operator
$$
\mathcal{O}:\,X(\Omega)\rightarrow
H^{s-2q,\,\varphi}_{\overline{\Omega}}(\mathbb{R}^{n}).
$$
It satisfies the condition $\mathcal{O}f=f$ in $\Omega$; i.e., $\mathcal{O}$ is
an operator extending a distribution from $\Omega$ onto $\mathbb{R}^{n}$. This
implies the continuous embedding $X(\Omega)\hookrightarrow
H^{s-2q,\,\varphi}(\Omega)$.

\begin{theorem} 
Let $s<2q-1/2$, $s+1/2\notin\mathbb{Z}$, and $\varphi\in\mathcal{M}$. We assume
that a Hilbert space $X(\Omega)$ is continuously embedded into
$\mathcal{D}'(\Omega)$ and satisfies Conditions $1$, $2$. Then the following
assertions hold true:
\begin{itemize}
\item[(i)] The set $D^{\infty}_{L,X}(\Omega):=
\{\,u\in C^{\infty}(\,\overline{\Omega}\,):\,Lu\in X(\Omega)\,\}$ is dense in
the space $D^{s,\varphi}_{L,X}(\Omega)$.
\item[(ii)] The mapping $(3.6)$, where $u\in D^{\infty}_{L,X}(\Omega)$,
is extended by a continuity to the linear bounded operator $(6.1)$.
\item[(iii)] The operator $(6.1)$ is a Fredholm one. Its kernel coincides with
$\mathcal{N}$, and its range is equal to the set
$$
\Bigl\{(f,g_{1},\ldots,g_{q})\in\mathcal{X}_{s,\varphi}(\Omega,\partial\Omega):\,
(f,v)_{\Omega}+\sum_{j=1}^{q}\,(g_{j},C^{+}_{j}v)_{\partial\Omega}=0\;\;\forall
\;\;v\in \mathcal{N}^{+}\Bigr\}.
$$
\item[(iv)] If the set $\mathcal{O}(X^{\infty}(\Omega))$ is dense in the space
$H^{s-2q,\,\varphi}_{\overline{\Omega}}(\mathbb{R}^{n})$, then the index of the
operator $(6.1)$ is equal to $\dim\mathcal{N}-\dim\mathcal{N} ^{+}$.
\end{itemize}
\end{theorem}

Conditions 1 and 2 allow us to vary the space $X(\Omega)$ in a broad fashion.
We especially note two possible options of $X(\Omega)$. The first of them is
the choice $X(\Omega):=H^{\sigma,\eta}(\Omega)$ for arbitrary fixed parameters
$\sigma>-1/2$ and $\eta\in\mathcal{M}$.

\begin{theorem} 
Let $s<2q-1/2$, $s+1/2\notin\mathbb{Z}$, $\sigma>-1/2$, and
$\varphi,\eta\in\mathcal{M}$. The mapping $(3.6)$ is extended by a continuity
to the bounded and the Fredholm operator
$$
(L,B):\,\bigl\{u\in H^{s,\varphi}(\Omega):Lu\in
H^{\sigma,\eta}(\Omega)\bigr\}\rightarrow H^{\sigma,\eta}(\Omega)
\oplus\bigoplus_{j=1}^{q}\,H^{s-m_{j}-1/2,\,\varphi}(\partial\Omega),
\eqno(6.2)
$$
provided that its domain is endowed with the graphics norm
$$
\bigl(\,\|u\|_{H^{s,\varphi}(\Omega)}^{2}+
\|Lu\|_{H^{\sigma,\eta}(\Omega)}^{2}\bigr)^{1/2}.
$$
The index of the operator $(6.2)$ is equal to
$\dim\mathcal{N}-\dim\mathcal{N}^{+}$ and does not depend on parameters $s$,
$\sigma$, $\varphi$, and $\eta$.
\end{theorem}

The case were $\sigma=0$ and $\eta\equiv1$, i.e.
$X(\Omega):=H^{0,1}(\Omega)=L_{2}(\Omega)$, is of great importance in the
spectral theory of elliptic operators \cite{Grubb68, Grubb96, Mikhailets82,
Mikhailets89}.

The condition $\sigma>-1/2$ is essential in Theorem 2, that does not allow us
to consider the boundary problem (3.1), (3.2) for an arbitrary distribution
$f\in \mathcal{D}'(\Omega)$ supported on a compact subset in $\Omega$. Here the
important example is $f(x):=\delta(x-x_{0})$, where $x_{0}\in\Omega$. The
following construction of the space $X(\Omega)$ has not this demerit.

We consider the set of weight functions
$$
\mathcal{W}^{\infty}_{k}(\,\overline{\Omega}\,):=\bigl\{\rho\in
C^{\infty}(\,\overline{\Omega}\,):\rho>0\;\;\mbox{in}\;\;\Omega,\;\;
D_{\nu}^{j}\,\rho=0\;\;\mbox{on}\;\;
\partial\Omega\;\;\forall\;\;j=0,\ldots,k\bigr\},
$$
where integer $k\geq0$.

Let $s<2q-1/2$, $\varphi\in\mathcal{M}$, and
$\rho\in\mathcal{W}^{\infty}_{[2q-s-1/2]}(\,\overline{\Omega}\,)$. (As usual,
$[t]$ denotes the integral part of $t$.) We consider the space
$$
\rho H^{s-2q,\,\varphi}(\Omega):=\bigl\{f=\rho v:\,v\in
H^{s-2q,\,\varphi}(\Omega)\,\bigr\}
$$
endowed with the inner product
$$
\bigl(f_{1},f_{2}\bigr)_{\rho H^{s-2q,\,\varphi}(\Omega)}:=
\bigl(\rho^{-1}f_{1},\rho^{-1}f_{2}\bigr)_{H^{s-2q,\,\varphi}(\Omega)}.
$$
The space $X(\Omega)=\rho H^{s-2q,\,\varphi}(\Omega)$ is Hilbert separable and
satisfies Conditions 1, 2.

\begin{theorem} 
Let $s<2q-1/2$, $s+1/2\notin\mathbb{Z}$, $\varphi\in\mathcal{M}$, and
$\rho\in\mathcal{W}^{\infty}_{[2q-s-1/2]}(\,\overline{\Omega}\,)$. The mapping
$(3.6)$, where $u\in C^{\infty}(\,\overline{\Omega}\,)$, $Lu\in\rho
H^{s-2q,\,\varphi}(\Omega)$, is extended by a continuity to the bounded and the
Fredholm operator
$$
(L,B):\,\bigl\{u\in H^{s,\varphi}(\Omega):Lu\in\rho
H^{s-2q,\,\varphi}(\Omega)\bigr\}\rightarrow \rho
H^{s-2q,\,\varphi}(\Omega)\oplus\bigoplus_{j=1}^{q}
\,H^{s-m_{j}-1/2,\,\varphi}(\partial\Omega), \eqno(6.3)
$$
provided that its domain is endowed with the graphics norm
$$
\bigl(\,\|u\|_{H^{s,\varphi}(\Omega)}^{2}+\|Lu\|_{\rho
H^{s-2q,\varphi}(\Omega)}^{2}\bigr)^{1/2}.
$$
The index of the operator $(6.3)$ is equal to
$\dim\mathcal{N}-\dim\mathcal{N}^{+}$ and does not depend on $s$, $\varphi$,
and $\rho$.
\end{theorem}

As an example of
$\rho\in\mathcal{W}^{\infty}_{[2q-s-1/2]}(\,\overline{\Omega}\,)$, we may chose
every function $\rho\in C^{\infty}(\,\overline{\Omega}\,)$ such that $\rho$ is
positive in $\Omega$ and
$$
\rho(\cdot)=(\mathrm{dist}(\cdot,\partial\Omega))^{\delta}\;\;\mbox{in a
neighborhood of}\;\;\partial\Omega\;\;\mbox{for}\;\;\delta=[2q-s+1/2].
\eqno(6.4)
$$

Theorems 6.1--6.3 were proved in \cite{MM15, MM16}. They are closely connected
with the theorems of J.-L.~Lions and E.~Magenes on a solvability of elliptic
boundary problems in the two-sided Sobolev scale \cite{LM62V, LM63VI, LM72,
Magenes65}. A theorem similar to Theorem 6.1 were proved in \cite[Sec.
6.10]{Magenes65} in the case of $s\leq0$, $\varphi\equiv1$ and the Dirichlet
boundary conditions. In this paper, certain different conditions depending on
the problem under consideration were imposed on $X(\Omega)$ (see also
\cite[Ch.~2, Sec. 6.2]{LM72}). Theorem 6.2 was proved in \cite{LM62V, LM63VI}
in the important case $\varphi\equiv\chi\equiv1$ and $\sigma=0$. Theorem 6.3
was proved in \cite[Ch. 2, Sec. 6,7]{LM72} in the case where $\varphi\equiv1$
and the weight function $\rho$ satisfies the condition (6.4) with
$\delta=2q-s$. The similar questions were considered in \cite{Roitberg68,
KostRoitb73}, \cite[Sec. 1.3]{Roitberg99} for the modified Sobolev scale. We
note that Theorems 6.2 and 6.3 are also true for half-integer values of $s$ if
we define the spaces with the help of the interpolation.

\end{document}